\documentclass[a4paper,12pt]{article}

\usepackage{latexsym,amssymb,amsmath}

\newtheorem{thm}{Theorem}
\newtheorem{lem}{Lemma}

\newtheorem{cor}{Corollary}

\newtheorem{rem}{Remark}
\newtheorem{exa}{Example}

\begin{document}

\title{A new integral formula for the inverse Fueter mapping theorem}

\author{Fabrizio Colombo$^{*}$\\
\small{e-mail: fabrizio.colombo@polimi.it}
\and Dixan Pe\~na Pe\~na$^{**}$\\
\small{e-mail: dpp@cage.ugent.be}
\and Irene Sabadini$^{*}$\\
\small{e-mail: irene.sabadini@polimi.it}
\and Frank Sommen$^{**}$\\
\small{e-mail: fs@cage.ugent.be}}

\date{\normalsize{$^*$Dipartimento di Matematica, Politecnico di Milano\\Via Bonardi, 9, 20133 Milano, Italy\\\vspace{0.2cm}
$^{**}$Clifford Research Group, Department of Mathematical Analysis\\Faculty of Engineering and Architecture, Ghent University\\Galglaan 2, 9000 Gent, Belgium}}

\maketitle

\begin{abstract}
\noindent In this paper we provide an alternative method to construct the Fueter primitive of an axial monogenic function of degree $k$, which is complementary to the one used in \cite{CoSaF2}. As a byproduct, we obtain an explicit description of the kernel of the Fueter mapping. We also apply our method to obtain the Fueter primitives of the Cauchy kernels with singularities on the unit sphere.\vspace{0.2cm}\vspace{0.1cm}
\end{abstract}

\section{Introduction}

Let us denote by $\mathbb{R}_{0,m}$ ($m\in\mathbb N$) the real Clifford algebra generated by the standard basis $\{e_1,\ldots,e_m\}$ of the Euclidean space $\mathbb R^m$  (see \cite{Cl}).
The multiplication in $\mathbb{R}_{0,m}$ is determined by the relations
\begin{equation*}
e_je_k+e_ke_j=-2\delta_{jk},\quad j,k=1,\dots,m
\end{equation*}
and a general element $a$ of $\mathbb R_{0,m}$ may be written as 
\[a=\sum_Aa_Ae_A,\quad a_A\in\mathbb R,\] 
in terms of the basis elements $e_A=e_{j_1}\dots e_{j_k}$, defined for every subset $A=\{j_1,\dots,j_k\}$ of $\{1,\dots,m\}$ with $j_1<\dots<j_k$. For the empty set, one puts $e_{\emptyset}=1$, the latter being the identity element. Conjugation in $\mathbb R_{0,m}$ is given by $\overline a=\sum_Aa_A\overline e_A$, where $\overline e_A=\overline e_{j_k}\cdots\overline e_{j_1}$ with $\overline e_j=-e_j$, $j=1,\dots,m$.

Observe that $\mathbb R^{m+1}$ may be naturally embedded in the real Clifford algebra $\mathbb R_{0,m}$ by associating to any element $(x_0,x_1,\ldots,x_m)\in\mathbb R^{m+1}$ the paravector given by 
\[x_0+\underline x=x_0+\sum_{j=1}^mx_je_j.\]
A function $f:\Omega\rightarrow\mathbb{R}_{0,m}$ defined and continuously differentiable in an open set $\Omega$ in $\mathbb R^{m+1}$ (resp. $\mathbb R^m$), is said to be monogenic if
\[(\partial_{x_0}+\partial_{\underline x})f=0\quad(\text{resp.}\;\partial_{\underline x}f=0)\;\;\text{in}\;\;\Omega,\]
where $\partial_{\underline x}=\sum_{j=1}^me_j\partial_{x_j}$ is the Dirac operator in $\mathbb R^m$ (see e.g. \cite{BDS,DSS,GuSp}). The differential operator $\partial_{x_0}+\partial_{\underline x}$, called generalized Cauchy-Riemann operator, gives a factorization of the Laplacian, i.e.
\[\Delta=\sum_{j=0}^m\partial_{x_j}^2=(\partial_{x_0}+\partial_{\underline x})(\partial_{x_0}-\partial_{\underline x}).\]
Thus the monogenic functions can be considered a subclass of the class of harmonic functions in $m+1$ variables.

Throughout the paper we assume  $h(z)=u(x,y)+iv(x,y)$ to be a holomorphic function in some open subset $\Xi$ of the upper half of the complex plane $\mathbb C$ and $P_k(\underline x)$ shall denote a homogeneous monogenic polynomial of degree $k$ in $\mathbb R^m$. Let us recall the following generalization of Fueter's theorem obtained in \cite{S}:
\begin{thm}
Put $\underline\omega=\underline x/r$, with $r=\vert\underline x\vert$, $\underline x\in\mathbb R^m$. If $m$ is odd, then the function
\begin{equation*}
\mathsf{Ft}\left[h(z),P_k(\underline x)\right](x_0,\underline x)=\Delta^{k+\frac{m-1}{2}}\bigl[\bigl(u(x_0,r)+\underline\omega\,v(x_0,r)\bigr)P_k(\underline x)\bigr]
\end{equation*}
is monogenic in $\Omega=\{(x_0,\underline x)\in\mathbb R^{m+1}:\;(x_0,r)\in\Xi\}$.
\end{thm}

In other words, this result provides a way to generate monogenic functions starting from a holomorphic function in the upper half of the complex plane. It was originally formulated by R. Fueter in the setting of quaternionic analysis in \cite{F} and later extended to the case of Clifford algebra-valued functions in \cite{Q,Sce} (see also  \cite{KQS,D,DS2,QS}).

\begin{rem}
{\rm It easily seen that $\mathsf{Ft}\left[h(z),P_k(\underline x)\right]$ defines an $\mathbb R$-linear operator between  holomorphic functions and monogenic functions considered
as real vector spaces, i.e.
\[\mathsf{Ft}\left[c_1h_1(z)+c_2h_2(z),P_k(\underline x)\right]=c_1\mathsf{Ft}\left[h_1(z),P_k(\underline x)\right]+c_2\mathsf{Ft}\left[h_2(z),P_k(\underline x)\right],\]
for all $c_1,c_2\in\mathbb R$.}
\end{rem}

The functions generated by this technique are monogenic functions of the form
\begin{equation}\label{AMF}
\bigl(A(x_0,r)+\underline\omega\,B(x_0,r)\bigr)P_k(\underline x),
\end{equation}
where $A$ and $B$ are $\mathbb R$-valued continuously differentiable functions in $\mathbb R^2$ which satisfy the following Vekua-type system
\begin{equation}\label{VESo}
\left\{\begin{array}{ll}\partial_{x_0}A-\partial_rB&=\displaystyle{\frac{2k+m-1}{r}}\,B\\\partial_{x_0}B+\partial_rA&=0.\end{array}\right.
\end{equation}
Monogenic functions of the form (\ref{AMF}) are called axial monogenic of degree $k$ and represent an important class of functions in Clifford analysis (see \cite{DSS}).

It is quite natural to ask whether given an axial monogenic function of degree $k$
\[H(x_0,\underline x)=\bigl(A(x_0,r)+\underline\omega\,B(x_0,r)\bigr)P_k(\underline x),\]
one can find a holomorphic function $h(z)$ such that
\[\mathsf{Ft}\left[h(z),P_k(\underline x)\right](x_0,\underline x)=H(x_0,\underline x).\]
The function $h$ is called the Fueter primitive of $H$. This problem has been recently studied in \cite{CoSaF,CoSaF2} and the Fueter primitive $h$ has been explicitly constructed. To this purpose, it was necessary to determine the Fueter primitives $\mathcal{W}^\pm _{k,m}$ of the functions
\[\mathcal{F}_{k,m}^+(x_0,\underline x)=\int_{\mathbb{S}^{m-1}}\mathcal{G}(x_0+\underline x-\underline{\omega})\,{P}_k(\underline{\omega})\,  dS(\underline{\omega}),\]
\[\mathcal{F}_{k,m}^-(x_0,\underline x)=\int_{\mathbb{S}^{m-1}}\mathcal{G}(x_0+\underline x-\underline{\omega})\, \underline{\omega}  \,\,{P}_k(\underline{\omega}) dS(\underline{\omega}),\]
where $dS(\underline{\omega})$ is the scalar element of surface area of $\mathbb{S}^{m-1}$,
and
\[\mathcal{G}(x_0+\underline x)=\frac{1}{A_{m+1}}\,\, \frac{\overline{x_0+\underline x}}{\vert x_0+\underline x\vert^{m+1}}\]
is the monogenic Cauchy kernel. Then, it is possible to express the Fueter primitive of $H$ in terms of a suitable integral involving $\mathcal{W}^\pm _{k,m}$, $A$, $B$ and $P_k$ (see \cite{CoSaF2}). This method can be used on any axially symmetric open set of $\mathbb R^{m+1}$, i.e. on every open set which is invariant under rotations that fix the real axis $x_0$.

Every monogenic function $f$ defined on an axially symmetric open set of $\mathbb R^{m+1}$ can be written as $f=\sum_{k=0}^\infty f_k$ where $f_k$ are axial monogenic functions of degree $k$. Hence for each term in the series we can provide a Fueter primitive as described above. We would like to note that the problem of inverting the Fueter mapping theorem has been recently tackled in the case of bi-axial monogenic functions (see \cite{CoSaF3}).

The aim of this paper is to present an alternative proof of the fact that the Fueter mapping is surjective on the set of axial monogenic functions of degree $k$ and to explicitly provide their Fueter primitives. The method we present here is complementary to the one presented in \cite{CoSaF2} in the sense that we here integrate with respect to the radius $r$ instead of the axial coordinate $x_0$. For the sake of simplicity the method is developed on a rectangle; of course it remains applicable on more general axially symmetric domains. 

As a byproduct of this method we describe the kernel of the Fueter mapping. We also compute with this method an exact formula for the Fueter primitives of the Cauchy kernels with singularities on the unit sphere. This corresponds to the integrals of the standard Cauchy kernel over the unit sphere.

\section{Some preliminary results}

Let $f:[a,b]\rightarrow\mathbb R$ be a continuous function. From the Cauchy formula for repeated integration, we know that an $n$-th antiderivative of $f$ is given by
\[f^{(-n)}(x)=\frac{1}{(n-1)!}\int_a^x(x-t)^{n-1}f(t)dt.\]
Inspired by this formula, we wish to find the solutions of the equations
\begin{equation}\label{Problem1}
\left(x^{-1}\frac{d}{dx}\right)^ng(x)=f(x)\quad\text{and}\quad\left(\frac{d}{dx}\,x^{-1}\right)^ng(x)=f(x),
\end{equation}
where $f:[a,b]\rightarrow\mathbb R$ is a given continuous function.

\begin{rem}
{\rm It is worth noting that the following identities hold (see \cite{HDPS,DS1})
\begin{equation*}
\left(x^{-1}\frac{d}{dx}\right)^ng(x)=\sum_{j=1}^n(-1)^{n+j}a_{j,n}x^{j-2n}\frac{d^jg}{dx^j}(x)
\end{equation*}
with
\[a_{j,n}=\frac{(2n-j-1)!}{2^{n-j}(n-j)!(j-1)!}\]
and
\begin{equation*}
\left(\frac{d}{dx}\,x^{-1}\right)^ng(x)=\sum_{j=0}^n(-1)^{n+j}a_{j+1,n+1}x^{j-2n}\frac{d^jg}{dx^j}(x).
\end{equation*}
Moreover, the integers
\[a_{j+1,n+1}=\frac{(2n-j)!}{2^{n-j}(n-j)!j!}\]
are the coefficients of the Bessel polynomial of degree $n$ (see \cite{KF}).}
\end{rem}

In order to find the solutions of (\ref{Problem1}), we define the following functions:
\begin{equation*}
\phi_n(x)=\int_a^xt\phi_{n-1}(t)dt,\quad\psi_n(x)=x\int_a^x\psi_{n-1}(t)dt,\quad x\in[a,b],\;n\ge1,
 \end{equation*}
with $\phi_0(x)=\psi_0(x)=f(x)$. Obviously, $\phi_n$ and $\psi_n$ satisfy
\begin{equation}\label{mainrr}
\frac{\phi_n^\prime(x)}{x}=\phi_{n-1}(x),\quad\left(\frac{\psi_n(x)}{x}\right)^\prime=\psi_{n-1}(x),\quad\phi_n(a)=\psi_n(a)=0,\;n\ge1.
\end{equation}

\begin{lem}\label{lemma1}
The functions $\phi_n$ and $\psi_n$ are given by
\begin{align}
\phi_n(x)&=\frac{1}{(2n-2)!!}\int_a^xt(x^2-t^2)^{n-1}f(t)dt,\label{Cfn1}\\
\psi_n(x)&=\frac{x}{(2n-2)!!}\int_a^x(x^2-t^2)^{n-1}f(t)dt,\label{Cfn2}
\end{align}
where $n!!$ denotes the double factorial of $n$.
\end{lem}
\textit{Proof.} Using integration by parts we obtain
\begin{multline*}
\int_a^xt(x^2-t^2)^{n-1}f(t)dt=\int_a^x(x^2-t^2)^{n-1}\phi_1^\prime(t)dt\\
=\big((x^2-t^2)^{n-1}\phi_1(t)\big)\big\vert_{t=a}^{t=x}+2(n-1)\int_a^xt(x^2-t^2)^{n-2}\phi_1(t)dt.
\end{multline*}
Then
\[\int_a^xt(x^2-t^2)^{n-1}f(t)dt=2(n-1)\int_a^xt(x^2-t^2)^{n-2}\phi_1(t)dt.\]
We iterate this procedure until the $(x^2-t^2)$ term vanishes. Thus after $n-1$ steps we
have
\[\int_a^xt(x^2-t^2)^{n-1}f(t)dt=2^{n-1}(n-1)!\int_a^xt\phi_{n-1}(t)dt=2^{n-1}(n-1)!\phi_n(x),\]
which proves (\ref{Cfn1}). Formula (\ref{Cfn2}) may be proved in a similar way.\hfill$\square$\vspace{0.2cm}

As an immediate consequence of Lemma \ref{lemma1} we obtain:
\begin{thm}\label{teoneces}
 Let $f:[a,b]\rightarrow\mathbb R$ be a continuous function.
The general solution of the equation
\[\left(x^{-1}\frac{d}{dx}\right)^ng(x)=f(x)\]
is
\[\phi_n(x)+\sum_{j=0}^{n-1}C_jx^{2j}\]
while the general solution of
\[\left(\frac{d}{dx}\,x^{-1}\right)^ng(x)=f(x)\]
is
\[\psi_n(x)+\sum_{j=0}^{n-1}\tilde{C}_jx^{2j+1},\]
where $C_j$ and $\tilde{C}_{j}$, $j=0,\dots,n-1$, are arbitrary real constants.
\end{thm}

\section{The inverse Fueter mapping theorem revisited}

Here plays an essential role the explicit form of $\mathsf{Ft}\left[h(z),P_k(\underline x)\right]$ determined in  \cite{D}, namely:
\begin{multline}\label{goodidea}
\mathsf{Ft}\left[h(z),P_k(\underline x)\right](x_0,\underline x)=(2k+m-1)!!\\
\times\left(\left(r^{-1}\partial_r\right)^{k+\frac{m-1}{2}}u(x_0,r)+\underline\omega\,\left(\partial_r\,r^{-1}\right)^{k+\frac{m-1}{2}}v(x_0,r)\right)P_k(\underline x).
\end{multline}

\begin{thm}[The inverse Fueter mapping theorem]
Let
\[H(x_0,\underline x)=\bigl(A(x_0,r)+\underline\omega\,B(x_0,r)\bigr)P_k(\underline x),\]
be a given arbitrary axial monogenic function of degree $k$ in
\[\Omega=\left\{(x_0,\underline x)\in\mathbb R^{m+1}:\;(x_0,r)\in [a,b]\times[c,d]\subset\mathbb R^2,\;c>0\right\}.\]
The Fueter primitives of $H(x_0,\underline x)$ exist and are given by
\begin{equation}\label{FPu}
u(x_0,r)=K_N\int_c^rt(r^2-t^2)^{N-1}A(x_0,t)dt+\sum_{j=0}^{N-1}\alpha_j(x_0)r^{2j},
\end{equation}
\begin{equation}\label{FPv}
v(x_0,r)=K_N\,r\int_c^r(r^2-t^2)^{N-1}B(x_0,t)dt+\sum_{j=0}^{N-1}\beta_j(x_0)r^{2j+1},
\end{equation}
where
\[K_N=\frac{1}{2N\left((2N-2)!!\right)^2}, \qquad N=k+\frac{m-1}{2}.\]
Moreover, the $\mathbb R$-valued functions  $\alpha_j(x_0)$ and $\beta_j(x_0)$ satisfy the following differential equations
\begin{multline}\label{cond1}
\alpha_j^\prime(x_0)-(2j+1)\beta_j(x_0)\\
=(-1)^{N-j-1}K_N\binom{N-1}{j}c^{2(N-j)-1}B(x_0,c),\;j=0,\dots,N-1,
\end{multline}
\begin{multline}\label{cond2}
\beta_j^\prime(x_0)+2(j+1)\alpha_{j+1}(x_0)\\
=(-1)^{N-j}K_N\binom{N-1}{j}c^{2(N-j-1)}A(x_0,c),\;j=0,\dots,N-2,
\end{multline}
\begin{equation}\label{cond3}
\beta_{N-1}^\prime(x_0)=-K_NA(x_0,c).
\end{equation}
\end{thm}
\textit{Proof.} Let $h(z)=u(x,y)+iv(x,y)$ be a Fueter primitive of $H$, i.e.
\[\mathsf{Ft}\left[h(z),P_k(\underline x)\right](x_0,\underline x)=H(x_0,\underline x),\]
then from (\ref{goodidea}) we obtain
\begin{align*}
\left(r^{-1}\partial_r\right)^{k+\frac{m-1}{2}}u(x_0,r)&=\frac{A(x_0,r)}{(2k+m-1)!!},\\
\left(\partial_r\,r^{-1}\right)^{k+\frac{m-1}{2}}v(x_0,r)&=\frac{B(x_0,r)}{(2k+m-1)!!}.
\end{align*}
Theorem \ref{teoneces} now yields (\ref{FPu}) and (\ref{FPv}). We must now investigate when $u$ and $v$ given by these formulae satisfy the Cauchy-Riemann equations, taking into account that  $A$ and $B$ fulfill the Vekua-type system (\ref{VESo}).

Let us define
\begin{equation}\label{Int1}
 I_1(x_0,r)=\int_c^rt(r^2-t^2)^{N-1}A(x_0,t)dt,
\end{equation}
\begin{equation}\label{Int2}
 I_2(x_0,r)=r\int_c^r(r^2-t^2)^{N-1}B(x_0,t)dt.
\end{equation}
It follows from (\ref{VESo}) and (\ref{mainrr}) that
\begin{align*}
\partial_{x_0}I_1(x_0,r)&=\int_c^rt(r^2-t^2)^{N-1}\left(\partial_tB(x_0,t)+\frac{2N}{t}B(x_0,t)\right)dt\\
&=-c(r^2-c^2)^{N-1}B(x_0,c)+(2N-1)\int_c^r(r^2-t^2)^{N-1}B(x_0,t)dt\\
&\qquad+2(N-1)\int_c^rt^2(r^2-t^2)^{N-2}B(x_0,t)dt,
\end{align*}
\begin{align*}
&\partial_{r}I_2(x_0,r)=\\
&\int_c^r(r^2-t^2)^{N-1}B(x_0,t)dt+2(N-1)r^2\int_c^r(r^2-t^2)^{N-2}B(x_0,t)dt=\\
&(2N-1)\int_c^r(r^2-t^2)^{N-1}B(x_0,t)dt+2(N-1)\int_c^rt^2(r^2-t^2)^{N-2}B(x_0,t)dt,
\end{align*}
which implies that
\[\partial_{x_0}I_1(x_0,r)-\partial_{r}I_2(x_0,r)=-c(r^2-c^2)^{N-1}B(x_0,c).\]
Similarly, we may verify that
\[\partial_{r}I_1(x_0,r)+\partial_{x_0}I_2(x_0,r)=r(r^2-c^2)^{N-1}A(x_0,c).\]
We thus have
\begin{multline*}
\partial_{x_0}u(x_0,r)-\partial_{r}v(x_0,r)=\sum_{j=0}^{N-1}\left(\alpha_j^\prime(x_0)-(2j+1)\beta_j(x_0)\right)r^{2j}\\
-K_N\,c\,(r^2-c^2)^{N-1}B(x_0,c),
\end{multline*}
\begin{multline*}
\partial_{r}u(x_0,r)+\partial_{x_0}v(x_0,r)=\sum_{j=0}^{N-2}\left(\beta_j^\prime(x_0)+2(j+1)\alpha_{j+1}(x_0)\right)r^{2j+1}\\
+\beta_{N-1}^\prime(x_0)r^{2N-1}+K_N\,r\,(r^2-c^2)^{N-1}A(x_0,c).
\end{multline*}
Therefore $u$ and $v$ satisfy the Cauchy-Riemann equations if and only if (\ref{cond1}), (\ref{cond2}) and (\ref{cond3}) are fulfilled.\hfill $\square$\vspace{0.2cm}

This theorem thus asserts that $\mathsf{Ft}\left[h(z),P_k(\underline x)\right]$ is surjective on the space of axial monogenic functions of degree $k$. Furthermore, we note that this operator is not injective since $\mathsf{Ft}\left[z^n,P_k(\underline x)\right](x_0,\underline x)=0$ for $0\le n\le 2k+m-2$, as it was observed in \cite{HDPS,DS1}.

In the next result, we show that the set of all real linear combinations of $z^n$, $0\le n\le 2k+m-2$, is indeed the kernel of $\mathsf{Ft}\left[h(z),P_k(\underline x)\right]$.

\begin{cor}
Let $\mathbb R_{2k+m-2}[z]$ be the vector space of all polynomials with real coefficients in $z$ of degree at most $2k+m-2$. Then
\[ker\left(\mathsf{Ft}\left[h(z),P_k(\underline x)\right]\right)=\mathbb R_{2k+m-2}[z].\]
\end{cor}
\textit{Proof.} We only have to prove that
\[ker\left(\mathsf{Ft}\left[h(z),P_k(\underline x)\right]\right)\subset\mathbb R_{2k+m-2}[z].\]
If $\mathsf{Ft}\left[h(z),P_k(\underline x)\right](x_0,\underline x)=0$, then from (\ref{FPu}) and (\ref{FPv}) we obtain
\begin{equation*}
u(x,y)=\sum_{j=0}^{N-1}\alpha_j(x)y^{2j},\quad v(x,y)=\sum_{j=0}^{N-1}\beta_j(x)y^{2j+1}.
\end{equation*}
The differential equations (\ref{cond1}), (\ref{cond2}) and (\ref{cond3}) now tell us that $\alpha_j(x)$ (resp. $\beta_j(x)$) are polynomials of degree at most $2(N-j)-1$ (resp. $2(N-j-1)$). Therefore
\[u(x,0)=C_0+C_1x+\dots C_{2N-1}x^{2N-1}, \quad v(x,0)=0,\]
for certain real constants  $C_0,\dots,C_{2N-1}$.  Then clearly $h(z)\in\mathbb R_{2k+m-2}[z]$.\hfill$\square$

\begin{rem}
{\rm As $\mathsf{Ft}\left[h(z),P_k(\underline x)\right]$  is an $\mathbb R$-linear operator, it is clear that
\[\mathsf{Ft}\left[h_1(z),P_k(\underline x)\right]=\mathsf{Ft}\left[h_2(z),P_k(\underline x)\right]\Leftrightarrow h_1(z)-h_2(z)\in\mathbb R_{2k+m-2}[z].\]}
\end{rem}
We end the paper with two examples.

\begin{exa}
{\rm From \cite{D,DS1} we have that
\[\mathsf{Ft}\left[1/z,P_k(\underline x)\right](x_0,\underline x)=(-1)^{k+\frac{m-1}{2}}((2k+m-1)!!)^2\left(\frac{x_0-\underline x}{\vert x_0+\underline x\vert^{2k+m+1}}\right)P_k(\underline x).\]
Thus if we apply (\ref{FPu}) and (\ref{FPv}) to this monogenic function we should be able to obtain the Cauchy kernel in the plane. Let us illustrate this for the case $k=0$, $m=5$. For this case $N=2$ and
\[A(x_0,r)=\frac{x_0}{(x_0^2+r^2)^3},\quad B(x_0,r)=-\frac{r}{(x_0^2+r^2)^3}.\]
It easily follows that
\[\beta_0(x_0)=-\frac{x_0^2+2c^2}{64(x_0^2+c^2)^2},\quad\alpha_0(x_0)=-x_0\beta_0(x_0),\]
\[\beta_1(x_0)=\frac{1}{64(x_0^2+c^2)^2},\quad\alpha_1(x_0)=-x_0\beta_1(x_0),\]
\[I_1(x_0,r)=\frac{(r^2-c^2)^2x_0}{4(x_0^2+r^2)(x_0^2+c^2)^2},\quad I_2(x_0,r)=-\frac{(r^2-c^2)^2r}{4(x_0^2+r^2)(x_0^2+c^2)^2},\]
where $I_1$ and $I_2$ denote the functions defined in (\ref{Int1}) and (\ref{Int2}), respectively.

Using (\ref{FPu}) and (\ref{FPv}) we obtain
\[u(x_0,r)=\frac{x_0}{64(x_0^2+r^2)},\quad v(x_0,r)=-\frac{r}{64(x_0^2+r^2)}.\]
That is $h(z)$ equals (up to a multiplicative constant) the Cauchy kernel in the plane.}
\end{exa}

\begin{exa}
{\rm In \cite{CoSaF} we have considered the functions
\[\mathcal{N}^+(q)=\int_{\mathbb{S}^{2}}\mathcal{G}(q-\underline{\omega})\,dS(\underline{\omega}),\]
and
\[\mathcal{N}^-(q)=\int_{\mathbb{S}^{2}}\mathcal{G}(q-\underline{\omega})\, \underline{\omega}\,dS(\underline{\omega}),\quad q=x_0+r\underline\omega,\]
and their Fueter primitives $\mathcal{W}^\pm$ in order to provide the Fueter inverse of a regular function of a quaternionic variable. This  corresponds to what we have discussed in the introduction, i.e. the integrals $\mathcal{F}_{k,m}^\pm$, in the particular case $k=0$, $m=3$.

In a closed form, these two functions can be written as
\[\mathcal{N}^+(q)=\frac{\underline\omega}{\pi r}\left(\frac{1}{1+q^2}-\frac{1}{r}{\rm Im}(\arctan q)\right),\]
and
\[\mathcal{N}^-(q)=\frac{\underline\omega}{\pi r}\left(\arctan q +\frac{q}{1+q^2}-\frac{1}{r}{\rm Im}(q\arctan q)\right).\]
We shall use formulae (\ref{FPu}) and (\ref{FPv}) to retrieve the fact that the Fueter primitive of $\mathcal{N}^+$ is the function $\mathcal{W}^+(z)=\frac{1}{2\pi}\arctan z$. Note that the function $\mathcal{N}^+$ is regular of axial type and it can thus be written as $\mathcal{N}^+(q)=A(x_0,r)+\underline{\omega}B(x_0,r)$ where
\[A(x_0,r)=\frac{1}{\pi}\frac{2x_0}{(1+x_0^2-r^2)^2+4x_0^2r^2}\]
and
\[B(x_0,r)=\frac{1}{2\pi r}\left(\frac{2\left(1+x_0^2-r^2\right)}{(1+x_0^2-r^2)^2+4x_0^2r^2}-\frac{1}{2r}\ln\left(\frac{x_0^2+(r+1)^2}{x_0^2+(r-1)^2}\right)\right).\]
We now compute $I_1$ and $I_2$ defined in (\ref{Int1}) and (\ref{Int2}) with $N=1$. We get that 
\[\begin{split}
2\pi I_1(x_0,r)&=2\pi\int_c^r tA(x_0,t)dt=\left.\arctan\left(\frac{2x_0}{1-x_0^2-t^2}\right)\right\vert_{t=c}^{t=r}\\
&=2{\rm Re}(\arctan z)-\arctan\left(\frac{2x_0}{1-x_0^2-c^2}\right)\\
&=2{\rm Re}(\arctan z) -4\pi\alpha_0(x_0),
\end{split}\]
\[\begin{split}
2\pi I_2(x_0,r)&= 2\pi r\int_c^rB(x_0,t)dt= \frac{r}{2t}\left.\ln\left(\frac{x_0^2+(t+1)^2}{x_0^2+(t-1)^2}\right)\right\vert_{t=c}^{t=r}\\
&=2{\rm Im}(\arctan z)-\frac{r}{ 2c}\ln\left(\frac{x_0^2+(c+1)^2}{x_0^2+(c-1)^2}\right)\\
&=2{\rm Im}(\arctan z)-4\pi r\beta_0(x_0).
\end{split}\]
On account of (\ref{FPu}) and (\ref{FPv}), we conclude that a Fueter primitive of $\mathcal N^+(q)$ is the function $\mathcal{W}^+(z)=\frac{1}{2\pi} \arctan z$ as it was computed in \cite{CoSaF}.

The function $\mathcal{N}^-(q)$ is also of axial type and so $\mathcal{N}^-(q)= A(x_0,r)+\underline{\omega}B(x_0,r)$ where
\[A(x_0,r)=\frac{1}{2\pi }\left(\frac{1}{2r} \ln\left(\frac{x_0^2+(r-1)^2}{x_0^2+(r+1)^2}\right)+\frac{2(x_0^2 +r^2-1)}{(1+x_0^2-r^2)^2+4x_0^2r^2}\right)\]
and
\[B(x_0,r)=\frac{x_0}{2\pi r}\left(\frac{1}{2r} \ln\left(\frac{x_0^2+(r-1)^2}{x_0^2+(r+1)^2}\right)+\frac{2(1+x_0^2+r^2)}{(1+x_0^2-r^2)^2+4x_0^2r^2}\right).\]
For this case we compute $I_1$ and $I_2$ with $N=1$ and obtain
\[\begin{split}
2\pi I_1(x_0,r)&=2\pi\int_c^r tA(x_0,t)dt\\
&=\left.\left(x_0\arctan\left(\frac{2x_0}{1-x_0^2-t^2}\right)-\frac{t}{2}\ln\left(\frac{x_0^2+(t+1)^2}{x_0^2+(t-1)^2}\right)\right)\right\vert_{t=c}^{t=r}\\
&=2 {\rm Re}(z\arctan z) -4\pi\alpha_0(x_0),
\end{split}\]
\[\begin{split}
2\pi I_2(x_0,r)&=2\pi r\int_c^rB(x_0,t)dt\\
&=r\left.\left(\frac{x_0}{2t}\ln\left(\frac{x_0^2+(t+1)^2}{x_0^2+(t-1)^2}\right)+\arctan\left(\frac{2x_0}{1-x_0^2-t^2}\right)\right)\right\vert_{t=c}^{t=r}\\
&=2{\rm Im}(z\arctan z)-4\pi r\beta_0(x_0).
\end{split}\]
Thus, using formulae (\ref{FPu}), (\ref{FPv})  we obtain that the Fueter primitive of $\mathcal{N}^-(q)$ is $\mathcal{W}^-(z)=\frac{1}{2\pi} z \arctan z$, as shown in \cite{CoSaF}.}
\end{exa}

\subsection*{Acknowledgments}

D. Pe\~na Pe\~na acknowledges the support of a Postdoctoral Fellowship funded by the \lq\lq Special Research Fund" (BOF) of Ghent University.

\end{document}